\title{CONWAY CLASSIFICATION OF ALTERNATING KNOTS}
\author{E. PI\~NA \\
Departamento de F\'\i sica\\ Universidad Aut\'onoma Metropolitana - Iztapalapa, \\
P. O. Box 55 534 \ Mexico, D. F., 09340 Mexico \\
e-mail: pge@xanum.uam.mx}

\documentclass[12pt]{article}
\usepackage{graphicx}
\usepackage{psfrag}
\begin{document}
\date {}
\maketitle

\abstract{The alternating knots, links and twists projected on the $S_2$ sphere were identified with the phase space of a Hamiltonian dynamic system of one degree of freedom. The saddles of the system correspond to the crossings, the edges correspond to the stable and unstable manifolds connecting the saddles. Each face is then oriented in one of two different senses determined by the direction of these manifolds. This correspondence can be also realized between the knot and the Poincar\'e section of a two degrees of freedom integrable dynamical system. The crossings corresponding to unstable orbits, and the faces to foliated torus, around a stable orbit.

The associated matrix to that connected graph was decomposed in two permutations. The separation was shown unique for knots not for links. The characteristic polynomial corresponding to some knot, link or twist families was explicitly computed in terms of Chebyschev polynomials. A classification of rational knots was formulated in terms of the first derivative of the polynomial of a knot computed in $x=2$, equal to the number of crossings of the knot multiplying the same number used previously by Conway for tabulation of knot properties. This leads to a classification of knots exemplified for the families having up to five ribbons. We subdivide the families of $N$ ribbons in subfamilies related to the prime knots of $N$ crossings.}

\newpage

\section{INTRODUCTION}
The importance of knots in the Dynamical Systems of Hamiltonian Dynamics is evident. Any periodic orbit of a Dynamical System in the Phase Space is the best example of a knot. If this example can be called applied or pure mathematics is a matter of taste. Considering the orbits of interacting particles in the presence of gravitational forces, one is conducted to consider periodic orbits forming a link of the relative positions of the closed orbits of several particles.

This study was restricted to the alternating knots that are defined when the knot is projected on a surface as a graph. One assumes overpasses and underpasses alternate along the curve. The restriction is extended to links with the same property independent of the number of closed curves. Most of the themes in this paper refers both to knots and links and we refer to them with the generic term of knot. When necessary we will distinguish the cases when we have exceptions using the term proper knot.

In Pi\~na (2006) one studied the matrix associated to the graph and its characteristic polynomial that are important semi-invariants of the knots.

The polynomials we are using are essentially different to several polynomials that were introduced in past in knot theory, like the Alexander's (1928), Jones' (1985), HOMFLY (see Freyd 1985), Kauffman's (1990),  Balister et Al. 2001, etc. These polynomials are recognized different to ours observing the polynomial corresponding to the unknotted circle has value zero whereas in other theories it is selected to have the value 1.

The classification used in this paper to identify particular knots follows the order presented by Rolfsen (1976) for prime knots and links, the first number denotes the number of crossings, the subindex is the place in the table of prime knots published in the Appendix C of the Rolfsen's book. This classification is the Standard one. For links the upper index denotes the number of components of the link in Appendix C of the same book.

When restricted to alternating knots it is not important what lane crosses in the upper level at any vertex. There are only two possibilities for the total knot, the one we choose and the reflected one, which can be identified with the mirror image of that knot and will not be distinguished in this paper. Once one lane is identified as upper at a crossing all the next crossings are determined by the alternating restriction.

We introduced an orientation of all the edges of the knot that is also found in Hamiltonian Dynamic Systems in the orientation of the stable and unstable manifolds connecting the saddles on a plane. The vertex is identified with the saddle and the upper strand with the unstable manifold, leaving the two edges the vertex with opposite orientation; the down strand is associated to the stable manifold formed by two edges that are oriented pointing to the vertex. This orientation of the edges give to any face one orientation clockwise or counterclockwise, allowing the construction of two Seifert's surfaces using one orientation to select the Seifert's circles and joining them at the crossing points by a bridge.
\begin{figure}

\psfrag{p}{\LARGE{$a_1$}}
\psfrag{q}{\LARGE{$a_2$}}
\psfrag{r}{\LARGE{$a_3$}}
\hfil\scalebox{0.5}{\includegraphics{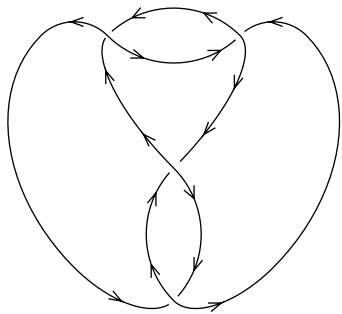}}\hfil

\caption{Knot of four crossings showing the orientation of the edges.}
\end{figure}

Balister et Al. (Balister 2001) have considered this orientation for alternating knots. In their terminology it is a connected two-in two-out digraph.

To see the graphs of a simple twist and the trefoil knot, oriented as explained, in a dynamic context see for example figure 3.2.5 at page 67 of Reichl (2004).

It is clear that the assumptions one adopted about these knots force some relative positions between knots and avoid to consider Reidemeister moves (Reidemeister 1948) of the edges of the knots.

\section{THE MATRIX OF THE ALTERNATING KNOTS AND THEIR POLYNOMIAL}
I studied in (Pi\~na (2006) the standard square matrix associated to an oriented graph (Berg\'e 2001). The matrix determines the alternating knot and vice versa. One enumerates the crossings. The edge connecting crossings from $j$ to $k$ is represented in the matrix by the entry $M_{jk} = 1$ at row $j$, column $k$ of the associated matrix. Exceptionally in some links two edges connect the same two crossings. In that cases the entry corresponding to those crossings will be equal to 2; the other entries are zero. The resulting matrix is an invariant of the knot, it can be thinked as equivalent to the knot. For Balister et Al. (Balister 2001) it is the {\sl adjacency} matrix. Generically we have two 1 and many 0 as entries of any row and any column of this matrix.

From this alternating matrix we identified if it is a knot or a link and how many different closed paths form it. We formed subgraphs starting at any non zero entry $M_{jk}$ of the associated matrix, join to the edge in the same row represented by the non zero entry $M_{jl}$, then we join the edge in the same column $M_{ml}$, then the edge in the same row $M_{mn}$ and so successively until the subgraph is closed with an edge $M_{zk}$ in the same column $k$ that the starting element. The resulting subgraph correspond to a knot or an unknotted element. A link formed by $K$ elements has $K$ subgraphs constructed in this way. Note that an entry 2 of the matrix correspond to a circle formed by two edges and appears only in links; by itself represent one component of the link.

The algebraic version (Pi\~na 2006) of the above properties is that the matrix of a knot is the sum of two permutation matrices, that were identified in alternated positions in the subgraph. This decomposition is unique for a proper knot. Each element of a link has a unique decomposition, but now the matrix is not decomposed in a unique form since the decompositions corresponding to each component in the link are independent.

Another property (Pi\~na 2006) was that any column and any row of the matrix has two 1's (or only one 2) and therefore the vector that has all the entries with the same value is an eigenvector of this matrix with eigenvalue 2. All the other eigenvectors are orthogonal to this vector, that means that the sum of all its components is null.

The characteristic polynomial of this matrix is a polynomial with integer coefficients since all the components of the matrix are integers. The roots of the polynomial are algebraic numbers. The polynomial factors always in two or more polynomials also having integer coefficients. Factor $(x - 2)$ appears in all the polynomials.

How good is this polynomial to identify a particular knot? The answer to this important question is not always positive. There are different links with same polynomial. And the same link with different positions of the components produces sometimes different polynomials. Nevertheless this polynomial is an important tool for deduce invariants of the alternating knots, in particular in this paper the first derivative of this polynomial at $x=2$ is an invariant related to the Alexander polynomial.

The matrices of many knots were studied and the characteristic polynomials of these matrices were computed and compared.

Many knots have an isolated face of two edges or a series of contiguous faces of two edges, forming a {\sl ribbon}. One notes that one face of two edges can be destroyed by eliminating the two edges and colliding or joining the two crossings in one vertex without modifying the rest of the knot and the orientation of the edges forming it. This operation can be iterated in a series of contiguous faces of two edges (ribbon), but also can be reversed to starting with a knot and obtaining other knot with an extra vertex in the same direction of the ribbon. This operation gives birth to a family of knots having the same structure but different number of faces of two edges in a contiguous series of a ribbon. These knots form a family and the characteristic polynomial of any member of the family can be computed as a function of three contiguous polynomials. Looking to the simplest families of polynomials one finds useful to follow the destruction of a series and also delete the last vertex by joining two times one incoming edge to an outgoing edge, the gap in the same direction of the ribbon. This is the so called (Farmer 1995) {\sl crossing elimination}, upon choosing from the two possibilities the ribbon direction. This produces in most cases other knot, but also produces a twist characterized by having an unknotted component with one or two loops, that could be unknotted. Also when we can separate two parts of a knot  by two edges, a {\sl composed} knot (Farmer 1995), we can twist these two edges and form a family of twist than could be unknotted by the characteristic Reidemeister move.

With this construction any vertex bifurcates in two families, with two possible orientations for the family. Or any vertex of a knot  without faces of two loops could be destroyed by thinking it as the last member of the family. These operations remind the Skein relation for the Kauffman polynomial (Kauffman 1990), other kind of polynomials associated to knots.

The matrix of a twist was computed in the same way we do for any knot. The only difference is that in some cases we have now loops: faces in a twist with the property of having only one edge. The polynomials of a family of twists are simpler that other family polynomials and we started with them.

We considered the simplest family of twist formed by twisting a circle. The first twist produces a figure eight with two loops. The next twist in this family has two loops joined to a face of two edges. The general term of the family is a twist with $V$ crossings, two loops, $V-1$ faces of two edges and the external face of $2V$ edges. The polynomials of this family obeys the recurrence relation that is found by induction
\begin{equation}
P_{V+2}(x) - x P_{V+1}(x) + P_V(x) = 0\, ,
\end{equation}
and can be determined by the two first. All these polynomials have the factor $x-2$ and are expressed as
\begin{equation}
P_V(x) = (x - 2) J_{V-1}(x)
\end{equation}
in terms of Polynomials $J_k(x)$ which are well known polynomials, particular cases of Jacobi, Gegenbauer and Chebyschev polynomials denoted as
\begin{equation}
J_{k}(x) = U_{k}(x/2)\, ,
\end{equation}
obeying the recurrence relation (1).
They can be written explicitly as
\begin{equation}
\left(\begin{array}{cc}
x & -1 \\
1 & 0
\end{array} \right)^{k+1} = \left(\begin{array}{cc}
J_{k+1}(x) & - J_k(x) \\
J_k(x) & - J_{k-1}(x)
\end{array} \right)\, .
\end{equation}
A different equivalent definition was presented in Pi\~na (2006). In some recurrences it is convenient to define $J_{-1}(x) = 0$.

I also studied (Pi\~na 2006) the twist families formed by twisting two edges that separate in two parts a knot. The knot with these two parts is the product of knots named composition (Farmer 1995) that give birth to the concept of {\sl prime knot} settled by H. Schubert in 1949. The polynomials of each family of twist satisfies the recurrence (1) that is a characteristic of the twists.

The polynomials of the families of knots obey a recurrence relation, that is different to the twist recurrence (1) by the fact that it has now a non-homogeneous term. But since all the polynomials have the factor $x-2$ it appears the same factor in the source, and the recurrence is always of the form
\begin{equation}
P_{V+2}(x) - x P_{V+1}(x) + P_V(x) = (x-2) H(x)\, ,
\end{equation}
where H(x) is a typical polynomial belonging to the family. This recurrence can be rewritten by transforming to the homogeneous form (1) by adding to any polynomial of the family the factor $H(x)$ of the source as
\begin{equation}
P_{V+2}(x) + H(x) - x [P_{V+1}(x) + H(x)] + [P_V(x) + H(x)] = 0 \, ,
\end{equation}
that can be solved as the twist examples. The polynomial $P_V(x)$ of a knot  is a linear combination of the Chebyshev polynomials minus the $H(x)$ polynomial.

We recover some examples (Pi\~na 2006). The family of cyclic torus knots begins with the eight shape twist of $V=1$, then the Hopf link, the trefoil knot, the Solomon's knot, $5_1$, $6_1^2$, $7_1$, $8_1^2$, etc. The polynomials of the knots  of this family were found to be of the form
\begin{equation}
2 [J_{V}(x) - 1] - x J_{V-1}(x)\, ,
\end{equation}
For $V$ an odd number we have a knot with polynomial
\begin{equation}
2 [J_{2 k + 1}(x)-1] - x J_{2 k} = (x - 2) [J_k(x) + J_{k - 1}(x)]^2\, .
\end{equation}
For $V$ an even number we have a link with
\begin{equation}
2 [J_{2 k}(x)-1] - x J_{2 k-1} = (x^2 - 4) [J_{k-1}(x)]^2\, .
\end{equation}

We found the polynomial of families formed by two ribbons $F_{j,k}(x)$ given by
$$
F_{j, k}(x) = (x-2) (J_{j-1}(x) + J_{k-1}(x))
$$
\begin{equation}
+ x \{ J_{j-1}(x) [J_k(x)-1] + [J_{j}(x)-1] J_{k-1}(x)\} -x^2 J_{j-1}(x) J_{k-1}(x) \, ,
\end{equation}

\begin{figure}

\hfil\scalebox{0.7}{\includegraphics{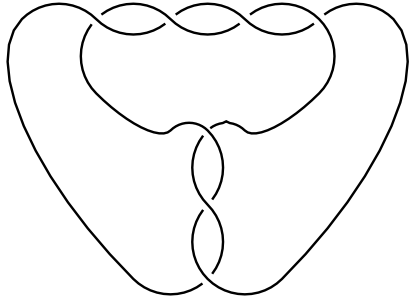}}\hfil

\caption{The knot $7_3$ formed by two ribbons.}
\end{figure}
\begin{figure}[b]

\psfrag{p}{\LARGE{4}}
\psfrag{q}{\LARGE{3}}
\hfil\scalebox{0.5}{\includegraphics{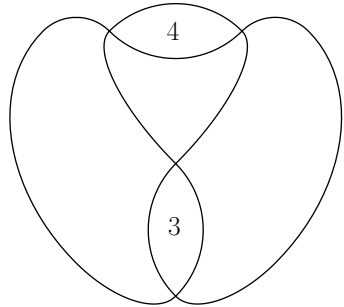}}\hfil

\caption{\small The same knot that in Fig. 2 with compact notation.}
\end{figure}

The non-symmetric three ribbon knot (see Fig. 4) produced the polynomial
$$
P_{k,l;m}(x) = \{J_k(x)J_l(x)+[x (J_{k-1}(x)-J_k(x)][x J_{l-1}(x) - J_{l}(x)]-2\}J_m(x)
$$
$$
- x \{[x J_{k-1}(x)-J_{k}(x)][x J_{l-1}(x)- J_{l}(x)]+ J_{k-1}(x)+ J_{l-1}(x)-1\}J_{m-1}(x)
$$
\begin{equation}%
- 2 J_{k-1}(x) J_{l-1}(x)\, ,
\end{equation}
symmetric in the $k$ and $l$ indexes, but non-symmetric in the $m$ index.
\begin{figure}

\psfrag{p}{\LARGE{$a_1$}}
\psfrag{q}{\LARGE{$a_2$}}
\psfrag{r}{\LARGE{$a_3$}}
\hfil\scalebox{0.5}{\includegraphics{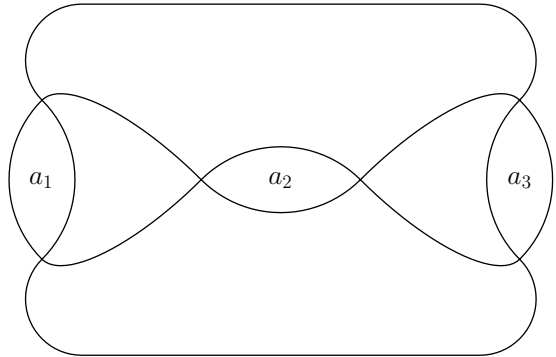}}\hfil

\caption{Knot formed by three ribbons two in parallel, one non parallel.}
\end{figure}
A similar knot formed by three ribbons (see Fig. 5) leads to the symmetric polynomial in the three indexes
$$
G_{k,l,m}(x) = x[J_{k-1}(x)J_l(x)J_m(x)+J_k(x)J_{l-1}(x)J_m(x)+J_k(x)J_l(x)J_{m-1}(x)]
$$
$$
-x^2[J_k(x)J_{l-1}(x) J_{m-1}(x)+J_{k-1}(x)J_l(x) J_{m-1}(x)+J_{k-1}(x) J_{l-1}(x)J_m(x)]
$$
\begin{equation}%
+ (x^3-2)J_{k-1}(x)J_{l-1}(x)J_{m-1}(x)-x[J_{k-1}(x)+J_{l-1}(x)+J_{m-1}(x)]\, .
\end{equation}
\begin{figure}

\psfrag{p}{\LARGE{$a_1$}}
\psfrag{q}{\LARGE{$a_2$}}
\psfrag{r}{\LARGE{$a_3$}}
\hfil\scalebox{0.5}{\includegraphics{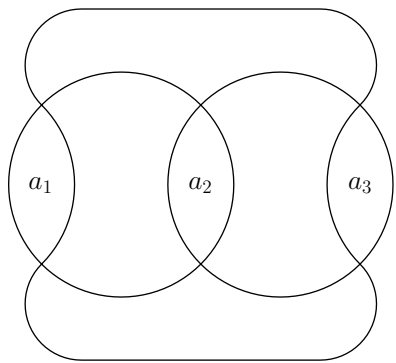}}\hfil

\caption{Knot formed by three ribbons in parallel.}
\end{figure}
Two rules are useful for allow the value zero in our equations. The uknotted and untwisted circle should have the polynomial with value zero. Two separated knots will have the polynomial equal to the product of the polynomials of the separated knots. This last rule includes the case when one or two have the polynomial equal to zero.

There are four different polynomials associated to two knots. The product, when they are unknotted. The composition when one edge of each is broken in two and reconnected between them forming one knot. One can form the first twist twisting the two edges of composition. A fourth case is when the two knots are linked forming a new two edges face with the same edges involved in the composition and the twist. The polynomials of the four cases are linearly related. Denoting the  polynomials of these four knots with the words: factor, composition, twist and link, the relation is
\begin{equation}%
link = (x+2) twist - (x+1) factor - (x) composition \, ,
\end{equation}
this equation and the previous recurrence between the members of a family of a ribbon produces many others polynomials, in particular the polynomials associated to composite knots have been computed. For example the composition of two and three ribbons in parallel gives the polynomials
\begin{equation}%
x [(J_k - 1) (J_{l-1} + J_{k-1} (J_l - 1)] - x^2 J_{k-1} J_{l -1} \, ,
\end{equation}
$$
(J_k - 1)(2 J_m - x J_{m-1} + 2)(J_l - 1)+ x^2 J_{k-1} J_{l -1}(J_m-xJ_{m-1})
$$
\begin{equation}%
- x [(J_k - 1) (J_{l-1} + J_{k-1} (J_l - 1)](J_m - x J_{m-1}+1) \, .
\end{equation}
Many other polynomials of families were computed and will be presented in the future if we discover a more efficient notation.

\section{THE CLASSIFICATION OF \\ ALTERNATING KNOTS}
A classification of knots is now formulated in terms of the number of different ribbons forming a knot and its relative connections. This classification is related closely to the study of {\sl tangles}, and classification of knots introduced by Conway (1970) in a classical paper on knots. What follows should be considered as a generalization or extension of the Conway's work.

The ribbons formed by 2 or more crossings determine a direction, the one joining the successive crossings that could be associated to the North-South or to the East-West direction. The two different directions are chosen according to the direction of the orientation of the external edges that connect the ribbon to the knot or the tangle such that the SE direction (also NW) is always outcoming, the SW and NE are incoming. If the ribbon has then the NS direction it is called vertical; otherwise it is horizontal. See in Figs. 1 and 2 one ribbon is horizontal, the other vertical.

Two ribbons are parallel if both are vertical or both horizontal. Two ribbons are orthogonal if one is vertical and the other horizontal.

Each ribbon is connected by four different edges to the rest of the knot. Cutting this portion of a knot gives an elementary tangle formed by one ribbon. Joining two of the edges of one elementary tangle formed by one ribbon of $k$ crossings to another elementary tangle formed by a ribbon of $l$ crossings produces a tangle formed by two ribbons in parallel or orthogonal directions, or an elementary tangle with one ribbon of $k + l$ crossings if the ribbons were connected one after the other. A general tangle is formed by an arbitrary number of ribbons. It has four edges disconnected.

Taking two ribbons with the same direction, they are connected one after the other if the resulting tangle is an elementary tangle formed by one ribbon with crossing number equal to the sum of the crossings of the two connected ribbons. In the next this connection is allowed only if it is considered as a one single ribbon tangle.

Two ribbons are connected in parallel if they are connected by two edges not one after the other (by the East or the West, if it is vertical; by the North or the South if it is horizontal), and appear both in a tangle with the same direction.

Connecting two edges of any tangle to two edges of another results in a more general tangle. Connection of tangles in alternating knots is allowed only respecting the direction of the edges, namely connecting the East of a left tangle with the West of a right tangle, or the North of a lower tangle with the South of an upper tangle.

Two knots (or two tangles) are equivalent (isotopic) by a flype. The successive connections by the East/West or North/South are obviously associative. This operation is also commutative between the connected tangles when all of them are of the East/West (or all of the North/South). In that case the commutation of two tangles is called a flype. Note that this property forces the union of different ribbons of the same direction when become one after the other by commutation.

Classification of alternating prime knots is by the number of different ribbons. An important helping function associated to a family is the first derivative of the corresponding polynomial computed in $x = 2$. Two properties were useful to translate from polynomial to derivative
\begin{equation}
\lim_{x\to 2} J_k(x) = k+1 \, , \quad \lim_{x\to 2} \frac {d J_k(x)}{d x}=\frac{k(k+1)(k+2)}{6}\, .
\end{equation}

The resulting derivative is a function of the literals representing the number of crossings in each ribbon. Since any polynomial has the factor $x - 2$ one obtains the same function dividing by $x - 2$ and taking the limit $x \rightarrow 2$. The resulting number is an invariant of the knot, equal to the number of crossings, multiplying an integer that is equal to the value of the Alexander polynomial at $x = -1$. This last number was tabulated by Conway (1970), it will be referred as Conway number (actually in this paper it is a positive number function of natural numbers,) and has the following properties.

\noindent The Conway function is a function of the number of crossings in each ribbon. It is a linear function on any of the crossing numbers. It is a polynomial with a number of terms, each of them are products of crossing numbers. The coefficient of any term is $+1$. In some cases the number $1$ is a term of a Conway number.

\noindent The number of twists do not appear in this number. Adding twists do not change the Conway number. The Conway number of a simple twist is 1. This is independent of the number of crossings in the twist.

\noindent The Conway number of the unknotted circle is defined one (notwithstanding its polynomial is zero).

\noindent The Conway number of the composition of two knots is equal to the product of the Conway numbers of the factors. Note the identity for composition is the unknotted circle.

\noindent The Conway number of two or more separated knots is zero. This results from factor $x - 2$ in each polynomial, or because the polynomial of the unknotted circle is zero.

\noindent The Conway number of the family of torus knots formed by one ribbon of $a_1$ crossings is $a_1$.

\noindent The Conway number of the family of prime knots formed by two orthogonal ribbons of crossings $a_1$ and $a_2$ is (see Figs. 2 and 3)
$$
1 + a_1 a_2 \, .
$$

We have two families of alternating prime knots with three ribbons. One (see Fig. 5) formed by three parallel ribbons of crossings $a_1$, $a_2$, $a_3$ with Conway number
$$
a_1 a_2 + a_2 a_3 + a_3 a_1 \, .
$$
The other (see Fig. 4) has two ribbons in parallel and one ribbon orthogonal to the other two with Conway number
$$
a_1 a_2 a_3 + a_1 + a_3 \, ,
$$
where $a_2$ is the crossing number of the different oriented ribbon. Compare these three simple expressions with the corresponding polynomials in equations (10-12).

A family of knots with $N$ ribbons can be sub-classified according to the set of knots of N crossings. The corresponding knot to a subfamily is obtained taking each crossing number equal to 1 in every ribbon of a family producing the representative knot. The Conway number of the representative knot is equal to the number of terms of the Conway number of the family. This results from taking equal to 1 the crossing numbers $a_j$. For example the two families of three ribbons correspond to the trefoil knot. The number of terms in both families is equal to the Conway number 3 of the trefoil knot.

The representative knot of $N$ crossings has $2 N$ edges connecting these crossings. These edges, their orientation, and relative position, will be present in any member of the associated family of $N$ ribbons that moreover could contain many other edges forming ribbons.
\begin{figure}

\psfrag{p}{\LARGE{$a_1$}}
\psfrag{q}{\LARGE{$a_2$}}
\psfrag{r}{\LARGE{$a_3$}}
\psfrag{s}{\LARGE{$a_4$}}
\hfil\scalebox{0.5}{\includegraphics{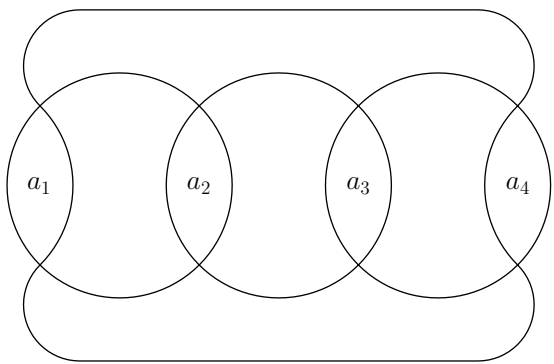}}\hfil

\caption{Knot formed by four ribbons in parallel.}
\end{figure}
The families of knots of 4 ribbons are 5, and can be separated according to the knot $4_1$ and the Solomon knot. The Solomon knot with Conway number 4 has two members, one (see Fig. 6) with 4 ribbons in parallel and Conway number
$$
a_1 a_2 a_3 + a_2 a_3 a_4 + a_3 a_4 a_1 + a_4 a_1 a_2 \, .
$$
\begin{figure}

\psfrag{p}{\LARGE{$a_1$}}
\psfrag{q}{\LARGE{$a_2$}}
\psfrag{r}{\LARGE{$a_3$}}
\psfrag{s}{\LARGE{$a_4$}}
\hfil\scalebox{0.5}{\includegraphics{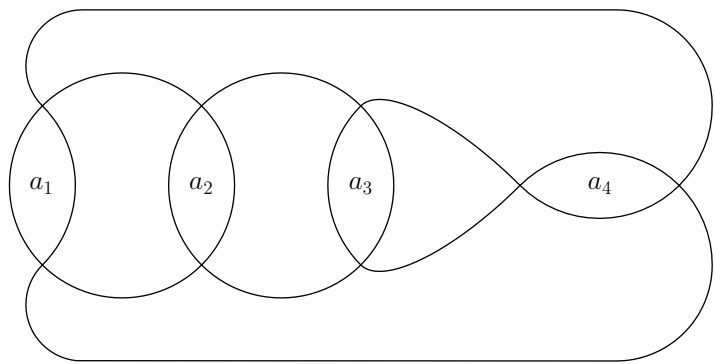}}\hfil

\caption{Knot formed by three ribbons in parallel, one non parallel.}
\end{figure}
The other (see Fig. 7) has 3 ribbons in parallel and one orthogonal to them. Its Conway number is
$$
a_1 a_2 a_3 a_4 + a_1 a_2 + a_2 a_3 + a_3 a_1 \, .
$$
We have not a family with 4 ribbons placed two parallel and two parallel orthogonal in alternating positions because that knot is equivalent to a member of the family of three ribbons with two parallel ribbons. This is a particular case of a more general property: any subfamily corresponding to a torus knot of 3 crossings or more, has only 2 families, the other possibilities are canceled by flypes.
\begin{figure}

\psfrag{p}{\LARGE{$a_1$}}
\psfrag{q}{\LARGE{$a_2$}}
\psfrag{r}{\LARGE{$a_3$}}
\psfrag{s}{\LARGE{$a_4$}}
\hfil\scalebox{0.5}{\includegraphics{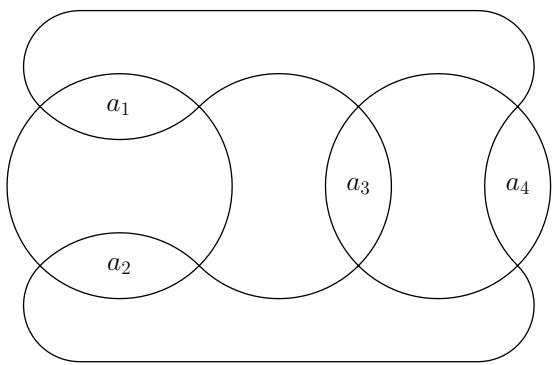}}\hfil

\caption{Knot formed by four ribbons in two parallel couples.}
\end{figure}
We have three families of 4 ribbons corresponding to the knot $4_1$. This knot has Conway number 5. First we connect a tangle formed by two parallel ribbons with crossing numbers $a_1$, $a_2$ to a tangle with two parallel ribbons orthogonal to the previous one with crossing numbers $a_3$, $a_4$. See Fig. 8. Its Conway number is
$$
a_1 a_2 a_3 a_4 + (a_1 + a_2)(a_3 + a_4) \, .
$$
\begin{figure}

\psfrag{p}{\LARGE{$a_1$}}
\psfrag{q}{\LARGE{$a_2$}}
\psfrag{r}{\LARGE{$a_3$}}
\psfrag{s}{\LARGE{$a_4$}}
\hfil\scalebox{0.5}{\includegraphics{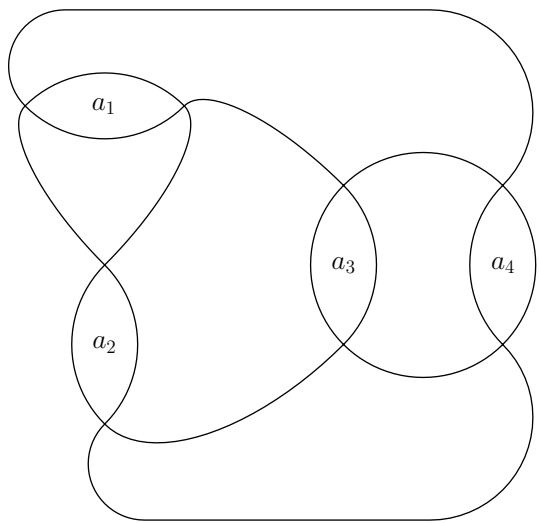}}\hfil

\caption{Knot formed by two parallel ribbons connected to a T tangle.}
\end{figure}
Second we connect one of the previous tangles with crossing numbers $a_1$, $a_2$ to the tangle formed by two orthogonal ribbons forming a $T$ tangle of crossing numbers $a_3$, orthogonal to the other three, and $a_4$. The relative positions of the two tangles should fit the relative positions of crossings of the $4_1$ knot. See Fig. 9. Its Conway number is
$$
a_1 a_2 a_3 + (a_1 + a_2)(a_3 a_4 + 1)\, .
$$
The third one family of 4 ribbons (See Fig. 10) is formed by two T tangles. A member of this family with 2 crossing numbers in each of the four ribbons is the knot $8_{12}$. The Conway number of this family is
$$
a_1 a_2 a_3 a_4 + a_1 a_2 + a_3 a_4 + a_1 a_4 + 1 \, .
$$
\begin{figure}

\

\psfrag{p}{\LARGE{$a_1$}}
\psfrag{q}{\LARGE{$a_2$}}
\psfrag{r}{\LARGE{$a_3$}}
\psfrag{s}{\LARGE{$a_4$}}
\hfil\scalebox{0.5}{\includegraphics{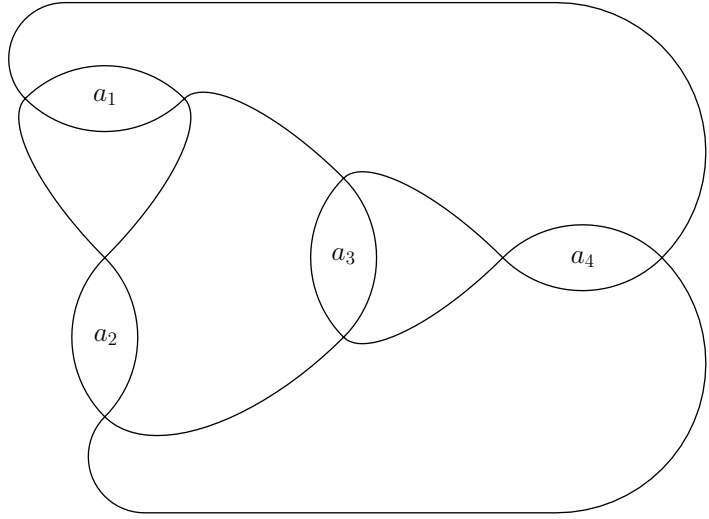}}\hfil

\caption{Rational knot formed by four ribbons.}
\end{figure}

The number of families of 5 ribbons is 12. Two of 5 terms associated to $5_1$, 4 of 7 terms associated to $5_2$ and 6 of 8 terms associated to the Whitehead link $5_1^2$. We list the Conway numbers of them.
$$
\begin{array}{l}
a_1 a_2 a_3 a_4 + a_2 a_3 a_4 a_5 + a_3 a_4 a_5 a_1 + a_4 a_5 a_1 a_2 + a_5 a_1 a_2 a_3 \\
a_1 a_2 a_3 + a_2 a_3 a_4 + a_3 a_4 a_1 + a_4 a_1 a_2 + a_1 a_2 a_3 a_4 a_5 \\
\ \\
(a_1 + a_2)(a_3 a_4 + a_4 a_5 + a_5 a_3) + a_1 a_2 a_3 a_4 a_5 \\
(a_1 a_2 + 1)(a_3 a_4 + a_4 a_5 + a_5 a_3) + a_2 a_3 a_4 a_5 \\
(a_1 + a_3 + a_1 a_2 a_3)(a_4 + a_5) + a_1 a_3 a_4 a_5 \\
(a_1 + a_3 + a_1 a_2 a_3)(1 + a_4 a_5) + a_1 a_3 a_5 \\
\ \\
1 + a_1 (a_2 + a_3) + (a_3 + a_4) a_5 + a_1 (a_2 a_3 + a_3 a_4 + a_4 a_2) a_5 \\
(a_1 + a_2)(a_4 + a_5) + a_1 a_2 a_3 (a_4 + a_5) + (a_1 + a_2)a_3 a_4 a_5 \\
(a_1 + a_2) a_3 (a_4 + a_5) + a_1 a_2 (a_4 + a_5) + (a_1 + a_2) a_4 a_5 \\
a_1 + a_5 + a_1 (a_2 + a_4) a_5 + (a_1 a_2 + 1) a_3 (a_4 a_5 + 1) \\
(a_1 + a_2)(a_3 a_4 a_5 + a_3 + a_5) + a_1 a_2 (1 + a_4 a_5) \\
(a_1 a_2 a_3 + a_1 + a_2)(1 + a_4 a_5) + (a_1 + a_2) a_4 a_5
\end{array}
$$

The number of factors in each term of the Conway number of a family is always an odd number or an even number. The difference is related to the parity of the maximum number of components in a link of the family. The proper knots having an even number of terms in the Conway number.

To construct these Conway functions one starts from the associated knot of $N$ crossings. One places $N$ ribbons with two possible orientations at each crossing; neglecting the cases when two ribbons are placed forming only one ribbon, the cases when a configuration is different only by different labeling of the ribbons, and the cases that are equivalent by a flype. After that, having a family, I make equal to zero successively the number of crossings of a particular ribbon. It comes a simplified family with known Conway number function. All the terms in this number also are terms of the number of the studied family. In some cases one obtain the Conway number by intersection of all these terms. An exception occurs when the product of all the crossing numbers $a_j$ is a term of the Conway number. That case is recognized since the total number of terms is equal to the Conway number of the associated knot of $N$ crossings.

Some times a Conway number is a function of two crossing numbers as the sum of them. That case correspond to a knot belonging to a family with lower number of ribbons. The knots are actually isotopic to the member of the family with lower number of ribbons, and they can be identified by using flypes and neglected.

One important particular case studied by many researchers (Kauffman 2004) is the set of rational knots that are in one-to-one correspondence with a continued fraction equal to the rate of two numbers: The Conway number is the numerator, corresponding to the rational knot, equal to the Gauss bracket of the $N$ crossing numbers $a_j$. The denominator is the Conway number of another rational knot, having $N-1$ ribbons, taking out the $a_1$ term, but keeping the same $a_j$ numbers for the rest of the ribbons, and the relative positions between them, that is characteristic of the rational knots. Actually (Kauffman 2004) reconnect the tangle so as to do the ribbon associated to the $a_1$ crossing to form a twist, resulting isotopic to a rational knot of crossing numbers $a_2$, $a_3$,..., $a_N$.

One notes that only one family of alternating knots of $N$ ribbons is rational. The knots with a small crossing number belong to several families when several crossing numbers of the ribbons are equal to 1. In such a case the ribbon has not a well defined orientation vertical or horizontal. This permits to select from various possible selections the one that is rational. When the number of crossing points of the knot is a larger number than 1 the exceptions to the rational cases become with plenty of possibilities. The five ribbons knots case is illustrative: from twelve families only one is rational. The cases associated to the $5_1$ and $5_2$ knots are never rational. If the knot is associated to the link $5_1^2$ is rational if it could be transformed to the rational form given at the last place in the forthcoming enumeration, otherwise it is not.

The first five families of rational knots have the Conway numbers equal to Gauss brackets as follows
$$
\begin{array}{l}
a_1 \\
1 + a_1 a_2 \\
a_1 a_2 a_3 + a_1 + a_3 \\
(a_1 a_2 + 1)(a_3 a_4 + 1) + a_1 a_4 \\
a_1 + a_5 + a_1 (a_2 + a_4) a_5 + (a_1 a_2 + 1) a_3 (a_4 a_5 + 1)
\end{array} \, .
$$

For 2, 3, 4 ribbons the rational families are represented by Figs. 3, 4, and 10, respectively.

One final point as a conclusion is to note that the Conway function presented in this paper could be defined by using the properties of this section, with no explicit relation to the adjacency matrix or the related polynomial trough the first derivative computed at $x=2$. For rational knots it coincides with the Gauss bracket. In general it is an important invariant for alternating knots that allow to distinguish and classify the families of alternating knots as are classified the rational knots by Gauss brackets.

\

{\large \bf ACKNOWLEDGMENTS }

\

This paper is dedicated to Prof. Donald Saari as a present. It should be presented at the Saarifest held in 2005 at Guanajuato, M\'exico. I thank to Prof. Lidia Jim\'enez-Lara for helping comments and outstanding literature.

\

{\large \bf REFERENCES }

\

\noindent C. C. Adams {\sl The Knot Book: An Elementary Introduction to the Mathematical Theory of Knots}, (American Mathematical Society, Providence, 2004)

\

\noindent J. W. Alexander "Topological Invariants of Knots and Links" {\sl Trans. Amer. Math. Soc.} {\bf 30} 275-306 (1928).

\

\noindent P. N. Balister, B. Bollob\'as, O. M. Riordan and A. D. Scott "Alternating Knot Diagrams, Euler Circuits and the Interlace Polynomial" {\sl Europ. J. Combinatorics} {\bf 22} 1-4 (2001).

\

\noindent C. Berg\'e {\sl The Theory of Graphs} (Dover, New York, 2001).

\

\noindent J. H. Conway "An enumeration of knots and links and some of their algebraic properties", in J. Leech (Ed.){\sl Computational Problems in Abstract Algebra}, (Pergamon Press, Oxford, 1970, pp. 329-358)

\

\noindent David. W. Farmer and T. B. Stanford {\sl Knots and Surfaces. A guide to Discovering} (American Mathematical Society, Providence, 1995).

\

\noindent P. Freyd, D. Yetter, J. Hoste, W. B. R. Likorish, K. Millet, and A. Ocneanu (HOMFLY) "A New Polynomial Invariant of Knots and Links" {\sl Bull. of the Am. Math. Soc.} {\bf 12} 239-246
(1985) .

\

\noindent V. F. R. Jones "A polynomial invariant for links via von Neumann algebra"{\sl Bull. of the Am. Math. Soc.} {\bf 12} 103-111 (1985).

\

\noindent L. H. Kauffman "An invariant of regular isotopy" {\sl Trans. Amer. Math. Soc.}, {\bf 318} 417-471 (1990).

\

\noindent L. H. Kauffman \& S. Lambropoulou "On the classification of rational tangles" {\sl Advances in Applied Mathematics} {\bf 33} (2004) 199-237.

\

\noindent E. Pi\~   na "Alternating knots and links theory"  arXiv:math.GT/060119v1 10 Jan 2006

\

\noindent L. E. Reichl {\sl The Transition to Chaos} 2nd ed. (Springer,
New York, 2004). p 67.

\

\noindent K. Reidemeister {\sl Knotentheorie} (Chelsea, New York, 1948).

\

\noindent D. Rolfsen {\sl Knots and Links}, (Publish or Perish, Berkeley, 1976).

\end{document}